\documentclass[a4paper,12pt]{article}
\usepackage[english]{babel}
\usepackage{amsmath}
\usepackage{amssymb}
\usepackage{tabls}
\usepackage[pdftex]{graphicx}
\usepackage{calc}
\usepackage{xcolor}
\usepackage[T1]{fontenc}
\usepackage[latin1]{inputenc}
\title{Generalized Maximal Orders}
\author{Tim Neijens\\University of Antwerp\\ \texttt{tim.neijens@gmail.com} \and Freddy Van Oystaeyen \\ University of Antwerp \\ \texttt{fred.vanoystaeyen@ua.ac.be}}

\newcommand{\Q}{\mathbb{Q}}

\newcommand{\blok}{\hfill \Box}
\newcommand{\CGR}{\mathop \diamondsuit \limits_{\sigma ,\alpha}}

\newtheorem{defi}{Definition}[section]
\newtheorem{opm}[defi]{Remark}
\newtheorem{prop}[defi]{Proposition}
\newtheorem{gev}[defi]{Corollary}
\newtheorem{st}[defi]{Theorem}
\newtheorem{lem}[defi]{Lemma}

\begin{document}

\maketitle
\begin{abstract}
Maximal Orders over an algebra are a generalization of the concept of a Dedekind domain.  The definition given in Maximal Orders by Reiner, \cite{R}, assumes that the field over which the algebra is defined is in the center of the order.  Since we want to define maximal orders over a Crystalline Graded Ring (defined in \cite{NVO6}), this concept needs to be generalized.  In this paper, we will weaken the condition that the field needs to be in the center, and still retain many of the desired properties of a maximal order.
\end{abstract}

\begin{flushleft}This paper is inspired by Maximal Orders, written by I. Reiner, \cite{R}.\end{flushleft}

\section{Preliminaries}

{\defi \label{def1}\textbf{Pre-Crystalline Graded Ring}\\
Let $A$ be an associative ring with unit $1_A$.  Let $G$ be an arbitrary group.  Consider an injection $u: G \rightarrow A$ with $u_e = 1_A$, where $e$ is the neutral element of $G$ and $u_g \neq 0$,  $\forall g \in G$.  Let $R \subset A$ be an associative ring with $1_{R}=1_A$.  We consider the following properties:
\begin{description} 
	\item[(C1)]\label{def2} $A = \bigoplus_{g \in G} R u_g$.
	\item[(C2)]\label{def3} $\forall g \in G$, $R u_g = u_g R$ and this is a free left $R$-module of rank $1$.
	\item[(C3)]\label{def4} The direct sum $A = \bigoplus_{g \in G} R u_g$ turns $A$ into a $G$-graded ring with $R = A_e$.
\end{description}
We call a ring $A$ fulfilling these properties a \textbf{Pre-Crystalline Graded Ring}.}\\

{\prop \label{def5} With conventions and notation as in Definition \ref{def1}:
\begin{enumerate}
	\item For every $g \in G$, there is a set map $\sigma_g : R \rightarrow R$   defined by: $u_g r = \sigma_g(r)u_g$ for $r \in R$.  The map $\sigma_g$ is in fact a surjective ring morphism.  Moreover, $\sigma_e = \textup{Id}_{R}$.
	\item There is a set map $\alpha : G \times G \rightarrow R$ defined by $u_g u_h = \alpha(g,h)u_{gh}$ for $g,h \in G$.  For any triple $g,h,t \in G$ the following equalities hold:
		\begin{eqnarray}
		\alpha(g,h)\alpha(gh,t)&=&\sigma_g(\alpha(h,t))\alpha(g,ht) \label{def6},\\
		\sigma_g(\sigma_h(r))\alpha(g,h)&=& \alpha(g,h)\sigma_{gh}(r) \label{def7}.
		\end{eqnarray}
	\item $\forall g \in G$ we have the equalities $\alpha(g,e) = \alpha(e,g) = 1$ and $\alpha(g,g^{-1}) = \sigma_g(\alpha(g^{-1},g)).$
\end{enumerate}
}
\begin{flushleft}\textbf{Proof}\end{flushleft} See \cite{NVO6}. $\blok$\\

{\prop Notation as above, the following are equivalent:
\begin{enumerate}
	\item $R$ is $S(G)$-torsionfree.
	\item $A$ is $S(G)$-torsionfree.
	\item $\alpha(g,g^{-1})r=0$ for some $g \in G$ implies $r = 0$.
	\item $\alpha(g,h)r=0$ for some $g,h \in G$ implies $r = 0$.
	\item $R u_g = u_g R$ is also free as a right $R$-module with basis $u_g$ for every $g \in G$.
	\item for every $g \in G$, $\sigma_g$ is bijective hence a ring automorphism of $R$.
\end{enumerate}
}
\begin{flushleft}\textbf{Proof}\end{flushleft} See \cite{NVO6}. $\blok$\\

{\defi Any $G$-graded ring $A$ with properties \textbf{(C1),(C2),(C3)}, and which is $G(S)$-torsionfree is called a \textbf{crystalline graded ring}.  In case $\alpha(g,h) \in Z(R)$, or equivalently $\sigma_{gh}=\sigma_g \sigma_h$, for all $g,h \in G$, then we say that $A$ is \textbf{centrally crystalline}.}\\

\section{General Theory}

\subsection{Setting} \label{max1} 
For this section, we will use the following notation and conventions.
\begin{itemize}
	\item $R$ is a Noetherian commutative domain.
	\item $K = \Q(R)$ the quotient field of $R$.
	\item $G$ is a finite group.
	\item There is a group action of $G$ on $R$, denoted by $\sigma_g \in \textup{Aut}R$, where $g \in G$.
	\item $A$ is a ring, and a vectorspace over $K$.  $K$ is not necessarily a subset of the center of $A$.
	\item $R^G \subset Z(A)$, where $R^G$ is the ring of invariant elements for the groupaction by $G$.  Extending for all $g\in G$ $\sigma_g$ to $K$, $K^G$ is the set of invariant elements in $K$ for the groupaction by $G$.
	\item For a set of generators $X_A$ of $A$ as a vectorspace over $K$, we have $Rx=xR, \forall x \in X_A$.  From this also follows $KX_A=A=X_AK$.
\end{itemize}

\subsection{Lattices}

{\defi We call $V$ a \textbf{$K$-bivectorspace over $K^G$} if $V$ is a left and right $K$-vectorspace and $arv = rva, \forall v \in V, a \in K^G, r \in K$.  We define an \textbf{$R$-bimodule over $R^G$} in an analoguous way.}\\

{\defi Let $V$ be a $K$-bivectorspace over $K^G$.  We call $M \subset V$ \textbf{an $R$-lattice} if $M$ is a left and right finitely generated $R$-torsionfree $R$-bimodule over $R^G$.  We call an $R$-lattice $M$ a \textbf{full $R$-lattice} if $KM=V=MK$.}\\

{\prop \label{max2} Let $M$ be an $R$-lattice, and $N$ be a full $R$-lattices in a $K$-bivectorspace over $K^G$ $V$.  Then $\exists r \neq 0 \in R^G$ such that $rM\subset N$, and $\exists s \neq 0 \in R^G$ such that $Ms \subset N$.}\\
\textbf{Proof} $N$ contains a left $K$-basis for $V$.  If $x \in M$, then we can choose $x = \sum a_i n_i$ for $n_i \in N$ and $a_i \in K$, $\forall i$.  We can choose $a_i$ so that the denominator is in $R^G$ (multiplying by the rest of the orbit of the $G$-action).  In other words $\exists r \in R^G$ such that $rx \in N$.  Since $M$ is finitely generated as an $R$-bimodule over $R^G$, we can choose $r \in R^G$ such that $rM \subset N$.\\
Since $N$ also contains a right $K$-basis for $V$, the second statement follows.$\blok$\\

\subsection{Orders} \label{max4}

{\defi We call a full $R$-latice $\Lambda$ in $A$, an \textbf{order in $A$} if it is a ring with the same identity element as $A$.  A full $R$-lattice that is maximal with this property is called a maximal order.}\\

{\defi \label{max22} Define for $M$, a full $R$-lattice in $A$ the \textbf{left order of $M$ $O_l(M)$} respectively the \textbf{right order of $M$ $O_r(M)$} as
\[O_l(M)=\{x \in A | xM \subset M\},\]
\[O_r(M)=\{x \in A | Mx \subset M\}.\]}

{\prop \label{max5} $O_l(M)$ and $O_r(M)$ as defined in Definition \ref{max22} are orders.}\\
\textbf{Proof} $O_l(M)$ and $O_r(M)$ are obviously rings, so we need to prove they are full $R$-lattices.\\
It is obvious that $\forall y \in X_A$, $yM$ is an $R$-lattice in $A$.  Indeed, it is an $R$-bimodule since $yR = Ry$.  So since $yM$ is an $R$-lattice and $M$ itself is full, we can use Proposition \ref{max2} to find an $r \in R^G$ such that $ryM \subset M$, or in other words $ry=yr \in O_l(M)$.  Since $X$ generates $A$ as a $K$-vectorspace, we find that $KO_l(M) =A=O_l(M)K$.\\
We still need to prove $O_l(M)$ is finitely generated.  $\exists s \in R$ with $s\cdot 1_D \in M$.  From this $O_l(M)\cdot(s \cdot 1_D)\subset O_l(M)M \subset M$.  So $O_l(M) \subset M s^{-1}$, and thus $O_l(M)$ is finitely generated.  Proving that $O_r(M)$ is a full $R$-lattice is analoguous, using the second statement of Proposition \ref{max2}. $\blok$\\

{\defi \label{max3} For any full $R$-lattice $M$ in $A$, define
\[M^{-1}=\{x \in A | MxM \subset M\}.\]}\\

{\prop $M^{-1}$ defined as above is a full $R$-lattice.}\\
\textbf{Proof} $M^{-1}$ obviously is an $R$-bimodule over $R^G$.  Clearly
\[L^{-1}=\{x \in A|Mx \subset O_l(M)\}=\{x \in A|xM\subset O_r(M)\}.\]
If $\Lambda = O_l(M)$, then $\Lambda$ is a full lattice since $M$ is (Proposition \ref{max5}).  Using Proposition \ref{max2} we find $\alpha, \beta \in R^G$ such that $\alpha \Lambda \subset M \subset \beta \Lambda$.  Therefore $\alpha^{-1} \Lambda \supset M \supset \beta^{-1} \Lambda$, since
\[\alpha^{-1}\Lambda = \{x \in A|\alpha \Lambda x \subset \Lambda\} \supset \{x \in A|Lx \subset \Lambda\}=L^{-1},\]
and likewise for the other inclusion.  Multiplying by $K$ gives us the desired result. $\blok$\\

For a full $R$-lattice $M_{ij}$ we define $\Lambda_i = O_l(M_{ij})$ and $\Lambda_j = O_r(M_{ij})$.

{\lem \label{max6} With notation as above, $M_{ij} \subset \Lambda_i \Leftrightarrow M_{ij} \subset \Lambda_j$}\\
\textbf{Proof} \[M_{ij} \subset \Lambda_i \Rightarrow M_{ij}M_{ij} \subset \Lambda_i M_{ij} = M_{ij} \Rightarrow M_{ij}\subset O_r(M_{ij})=\Lambda_j.\]
The other arrow is similar. $\blok$\\

{\lem \label{max7}Let $L$ be a full $R$-lattice.  Then
\[O_r(L^{-1})\supset O_l(L), \ \ \ O_l(L^{-1}) \supset O_r(L).\]}
\textbf{Proof} Let $x \in L^{-1}$.  Then
\begin{eqnarray*}
LxL \subset L &\Rightarrow& Lx \subset O_l(L)\\
\ &\Rightarrow& LxO_l(L) \subset O_l(L)\\
\ &\Rightarrow& LxO_l(L)L\subset O_l(L)L \subset L\\
\ &\Rightarrow& xO_l(L)\subset L^{-1}.
\end{eqnarray*}
This implies that $L^{-1}$ is a right $O_l(L)$-module.  So $L^{-1}O_l(L) \subset L^{-1}$ and thus $O_l(L) \subset O_r(L^{-1})$.  The other inclusion is similar. $\blok$\\

\subsection{Ideals}

{\opm All ideals $M$ defined below are assumed to be full $R$-lattices in $A$ and written "Ideal".}\\

{\defi Call $M$ a \textbf{normal Ideal} if $O_l(M)$ is a maximal $R$-order.  An \textbf{integral Ideal} is a normal Ideal such that $M \subset O_l(M)$.  A \textbf{maximal integral Ideal} is an integral Ideal $M$ which is a maximal left Ideal in $O_l(M)$.}\\

{\defi The normal Ideal $M$ is \textbf{two-sided} if $O_l(M)=O_r(M).$}\\

We are now able to repeat the classical theory of maximal orders, in an extended set-up.  The proofs are close or the same to those in \cite{R}.\\

{\defi\label{max8} A \textbf{prime Ideal} of an $R$-order $\Lambda$ is a proper two-sided Ideal $P$ in $\Lambda$ such that for each pair of two-sided Ideals $S$ and $T$ in $\Lambda$ with $ST\subset P$ we have $S \subset P$ or $T \subset P$.}\\

{\prop \label{max9} With notations as in Definition \ref{max8}, we can assume both $S$ and $T$ contain $P$ and as such $J \cdot J'=0$ in $\Lambda/P$, then $J=0$ or $J'=0$.}\\
\textbf{Proof} Elementary. $\blok$\\

{\st \label{max10} Let $R$ be a Dedekind domain.  The prime Ideals in $\Lambda$ coincide with the maximal two-sided Ideals of $\Lambda$.  If $P$ is a prime Ideal of $\Lambda$ then $p = P \cap R^G$ is a nonzero prime in $R^G$ and $\bar \Lambda = \Lambda/P$ is a finite dimensional simple $R^G/p$-ring.}\\
\textbf{Proof} Let $P$ be a prime Ideal in $\Lambda$ and set $p = R^G \cap P$, $\bar \Lambda = \Lambda/P$.  $P$ is a full $R$-lattice in $\Lambda$, so $p \neq 0$.  $p$ is a proper Ideal of $R^G$ since $1 \notin P$.  Now if $\alpha, \beta \in R^G$ then

\begin{eqnarray*}
\alpha\beta \in p &\Rightarrow& \alpha\Lambda\beta\Lambda \subset P \Rightarrow \alpha\Lambda \subset P \textup{ or } \beta\Lambda \subset P\\
\ &\Rightarrow& \alpha \in p \textup{ or } \beta \in p,
\end{eqnarray*}
so $p$ is prime and $\bar \Lambda$ is finite dimensional over $R^G/p$.  Since $\bar \Lambda$ is Artinian, $\textup{rad}\bar\Lambda$ is nilpotent and thus equal to $0$ by Proposition \ref{max9}.  $\bar \Lambda$ is semisimple.  The simple components are two-sided Ideals that annihilate each other, so again by Proposition \ref{max9} there is only one simple component.  This means $\bar \Lambda$ is simple and $P$ is maximal.  The other way is trivial. $\blok$\\

{\lem \label{max11} Let $M$ be an Ideal in a maximal order $\Lambda$ of $A$.  Then $O_l(M)=O_r(M)=\Lambda$.  This implies that $M^{-1}$ is a $\Lambda$-Ideal in $A$.  Furthermore, if $N$ is another Ideal in $\Lambda$ with $N \subset M$, then $M^{-1}\subset N^{-1}$.}\\
\textbf{Proof} Remark that $\Lambda \subset O_l(M)$ and $\Lambda \subset O_r(M)$, and $\Lambda$ is a maximal order.  We have from Lemma \ref{max7} that $O_l(M^{-1}) \supset \Lambda$ and $O_r(M^{-1}) \supset \Lambda$.  Suppose $N \subset M$, and take $x \in M^{-1}$.  Then $MxM \subset M$ and thus $xM \subset O_r(M)=\Lambda$.  So $NxN \subset NxM \subset N\Lambda \subset N$.  $\blok$\\

{\opm With notation as in Lemma \ref{max11}, we find
\[M^{-1} = \{x \in A| Mx \subset \Lambda\}=\{x \in A | xM \subset \Lambda\}.\]}\\

{\lem \label{max12} Let $M$ be a two-sided Ideal in the maximal order $\Lambda$.  If $M \neq \Lambda$, then $M^{-1} \neq \Lambda$.}\\
\textbf{Proof} Suppose that $M^{-1}=\Lambda$.  Since $M \neq \Lambda$, $\exists P$ prime in $\Lambda$ such that $M \subset P \subset \Lambda$.  Let $\alpha$ be a nonzero element in $R^G \cap P$.  (This exists, since $P$ is a full $R$-lattice (it contains $M$) in $A$ and we can multiply each element in $R$ with the other elements in the orbit under action of $G$ to obtain an invariant element.)  Since we are working in Noetherian rings, the two-sided Ideal $\alpha \Lambda$ contains a product of prime ideals:
\[P \supset \alpha\Lambda \supset P_1 \cdot \ldots \cdot P_r,\]
where the $P_i$ are prime ideals and $r$ is minimal.  Since $P$ is prime, it follows that $P=P_j$ for some $j$.  We may then write
\[P \supset \alpha\Lambda \supset BPC,\]
where $B,C$ are two-sided Ideals in $\Lambda$.  Thus
\begin{eqnarray*}
\alpha^{-1}BPC \subset \Lambda &\Rightarrow& B \alpha^{-1}PC B \subset B\\
\ &\Rightarrow& \alpha^{-1}PCB \subset O_r(B)=\Lambda\\
\ &\Rightarrow& P\alpha^{-1}CBP \subset P\\
\ &\Rightarrow& \alpha^{-1}CB \subset P^{-1} \subset M^{-1} = \Lambda.
\end{eqnarray*}
And this shows that $\alpha \Lambda$ contains a product of $r-1$ primes, contradiction. $\blok$\\

{\st \label{max13} Let $M$ be a full $R$-lattice in $A$ such that $O_l(M)$ is a maximal order.  Then
\[M \cdot M^{-1} = O_l(M).\]}\\
\begin{flushleft}\textbf{Proof}\end{flushleft}  Let $\Lambda=O_l(M)$ be maximal, set $B = MM^{-1}$.  From Lemma \ref{max11} we find that $B$ is a two-sided Ideal in $\Lambda$.  And so
\[BB^{-1} \subset \Lambda \Rightarrow MM^{-1}B^{-1} \subset \Lambda \Rightarrow M^{-1}B^{-1} \subset M^{-1}\rightarrow B^{-1} \subset O_r(M^{-1}).\]
But $O_r(M^{-1})\supset O_l(M) = \Lambda$ by Lemma \ref{max7}.  So from maximality, $\Lambda = O_r(M^{-1})$.  Therefore $B^{-1}\subset \Lambda$ and now $B = \Lambda$ from Lemma \ref{max12}.  $\blok$\\

{\st \label{max14} Let $\Lambda$ be a maximal order.  For each two-sided Ideal $M$ in $\Lambda$ we have
\[MM^{-1} = M^{-1}M = \Lambda,\]
\[(M^{-1})^{-1}=M.\]
Furthermore, every two-sided Ideal in $\Lambda$ is uniquely expressible as a product of prime Ideals of $\Lambda$, and multiplication of prime Ideals is commutative.}\\
\textbf{Proof}  The first formula is clear from Theorem \ref{max13}.  For the second formula, choose a nonzero $\alpha \in R^G$ such that $\alpha M^{-1} \subset \Lambda$.  This is a two-sided Ideal of $\Lambda$, so we have
\[(\alpha M^{-1})(\alpha M^{-1})^{-1}=\Lambda.\]
It is easily checked that
\[(\alpha M^{-1})^{-1} = \alpha^{-1}(M^{-1})^{-1}.\]
This implies
\[M^{-1}L = \Lambda \Rightarrow (M^{-1})^{-1}= \Lambda (M^{-1})^{-1} = M M^{-1}(M^{-1})^{-1}=M \Lambda = M.\]
We shall now prove that every two-sided Ideal in $\Lambda$ is a product of prime Ideals, or an empty product, in which case the Ideal is equal to $\Lambda$.  Suppose this is false, and let $N$ be a maximal counterexample, this itself is not prime.  So there exists a prime Ideal $P$ with $N \subset P \subset \Lambda$ and $N \neq P$.  So
\[N = N \Lambda \subset NP^{-1} \subset NN^{-1} = \Lambda.\]
If $N = NP^{-1}$ then $P^{-1} \subset O_r(N)=\Lambda$ and this is impossible by Lemma \ref{max12}.  So $NP^{-1}=P_1\cdot \ldots \cdot P_r$.  This implies $N = P_1\cdot \ldots \cdot P_r P$ which is impossible.  So every two-sided Ideal is a product of prime Ideals.\\
Now we want to establish commutativity.  Let $P, P'$ be distinct prime Ideals of $\Lambda$. Then
\[P^{-1}P'P \subset P^{-1}\Lambda P = \Lambda\]
and
\[P(P^{-1}P'P)=P'P \subset P'.\]
Thus $P'$ contains the product $P(P^{-1}P'P)$ of a pair of two-sided Ideals of $\Lambda$ and primality gives us $P' \supset (P^{-1}P'P)$, or $PP' \supset P'P$.  Since the reverse inclusion follows from symmetry, we obtain $PP' = P'P$.\\
We now prove the uniqueness in decomposition (up to the order of the factors).  So let
\[P_1 \cdot \ldots \cdot P_r = Q_1 \cdot \ldots \cdot Q_s,\]
where the $P_i$ and $Q_j$ are prime Ideals.  Then $P_1 \supset \prod Q_j$ implies $P_1 = Q_k$ for some $k$.  We then multiply both sides by $(P_1)^{-1}$ and repeat the process.  $\blok$\\

From this we find that the set of two-sided $\Lambda$-Ideals in $A$ is the free abelian group generated by the prime Ideals of $\Lambda$.\\

From now on, we suppose $R$ is a Dedekind domain.\\

{\gev \label{max15} If $M$ is a proper left Ideal of the maximal order $\Lambda$, then $M^{-1} \neq \Lambda$.}\\
\textbf{Proof} Let $N$ be a maximal left Ideal of $\Lambda$, so we find $M \subset N \subset \Lambda$.  Modifying the proof of Lemma \ref{max11} we find that $M^{-1} \supset N^{-1} \supset \Lambda$.  It suffices to prove that $N^{-1} \neq \Lambda$.  Set
\[P = \textup{ann}_\Lambda \Lambda/N = \{x \in \Lambda | x\Lambda \subset N\}.\]
This is a two-sided Ideal of $\Lambda$.  Choose a nonzero $\alpha \in R^G$ such that $\alpha \Lambda \subset N$. then $\alpha \in P$, so we find $\alpha \Lambda \subset P$.  So we see that $\Lambda/P$ is an artinian ring.  Now $\Lambda/N$ is a faithful left $(\Lambda/P)$-module, and is simple as such.  Therefore the ring $\bar\Lambda =\Lambda/P$ is a simple artinian ring and this makes $P$ into a prime Ideal.  Hence $N$ is the inverse image under the map $\Lambda \rightarrow \bar\Lambda$ of a maximal left Ideal of $\bar \Lambda$.  Such an Ideal is of the form
\[\bar \Lambda (1-f),\]
where f is a primitive idempotent in $\bar \Lambda$.  So
\[N = \Lambda(1-e) + P,\]
for some $e \in \Lambda$ which maps to $f$.  Set $L = e\Lambda + P$, a right Ideal in $\Lambda$.  Since $(1-e)e \in P$, $NL \subset P$.  $NL$ also is a two-sided Ideal that contains $P^2$, so $NL = P^2$ or $NL =P$.  If $NL = P^2$ then $Pe \subset P^2$ and multiplying by $P^{-1}$ we see that $e \in P$, contradiction.  This shows that $NL = P$.\\
Choose a nonzero $\alpha \in R^G \cap P$, and write $\alpha \Lambda$ as a product of prime Ideals of $\Lambda$.  Since $\alpha \Lambda \subset P$, one of the factors must be equal to $P$.  Thus we can write, with $Q$ some two-sided Ideal of $\Lambda$
\[\alpha\Lambda = PQ = NLQ.\]
Suppose $N^{-1} = \Lambda$ then
\[\alpha^{-1}LQ \subset N^{-1} = \Lambda \Rightarrow LQ \subset \alpha\Lambda \Rightarrow L\Lambda \subset \alpha\Lambda Q^{-1}=P \Rightarrow L \subset P.\]
But that last inclusion is impossible, so $N^{-1} \neq \Lambda$. $\blok$\\

{\st Let $M$ be a full $R$-lattice in $A$.  Then $O_l(M)$ is a maximal order if and only if $O_r(M)$ is a maximal order.}\\
\textbf{Proof} Let $\Lambda = O_l(M)$ be maximal.  Replacing $M$ by $\alpha M$ with $\alpha \in R^G$ and nonzero, we may assume that $M \subset \Lambda$.  We can assume that $M \neq \Lambda$, otherwise the result is trivial.  Suppose that the theorem is false, take $M$ a maximal counterexample.  Choose a left Ideal $L$ of $\Lambda$ such that $M \subset L \subset \Lambda$, $M \neq L$ and $L/M$ is a simple left $\Lambda$-module.  Then $O_l(L)=\Lambda$ and since $M$ was a maximal counterexample, $O_r(L) = \Lambda'$ is maximal order.  By Lemma \ref{max13} $LL^{-1} = \Lambda$.  By the same Lemma, with a little modification, we find $L^{-1}L = \Lambda'$.\\
Each left Ideal $X$ of $\Lambda'$ maps onto a $\Lambda$-submodule $LX$ of $L$, and the correspondence $X \leftrightarrow LX$ is one-to-one, since $X = L^{-1}LX$.  Set $N = L^{-1}M$, then
\[N = L^{-1}M \subset L^{-1}L = \Lambda'.\]
Therefore $N$ is a maximal left Ideal of $\Lambda'$ because $M$ is a maximal $\Lambda$-submodule of $L$.  By Lemma \ref{max13} and Corollary \ref{max15} we obtain
\[N^{-1} \supset \Lambda' ,\ \ \ N^{-1} \neq \Lambda',\]
\[N N^{-1}=\Lambda',\]
\[N^{-1}N \subset O_r(N).\]
Let $x \in A$ then
\[Nx \subset N \Rightarrow LNx \subset LN \Rightarrow x \in O_r(M),\]
\[Mx \subset M \Rightarrow L^{-1}Mx \subset L^{-1}M \Rightarrow x \in O_r(N).\]
This implies $O_r(M)=O_r(N)$, and from the choice of $M$, this is not a maximal order.  Thus there exists a maximal order $\Lambda''$ with $O_r(N) \subset \Lambda''$ and $O_r(N) \neq \Lambda''$.  Since
\[N \Lambda'' N^{-1} \cdot N\Lambda''N^{-1} \subset N \Lambda''O_r(N)\Lambda''N^{-1}\subset N \Lambda'' N^{-1},\]
$N \Lambda'' N^{-1}$ is an order.  Also
\[N \Lambda'' N^{-1} \supset NN^{-1}=\Lambda',\]
$\Lambda'$ is maximal so $N \Lambda'' N^{-1} = \Lambda'$.  So we have the following chain of $\Lambda'$ Ideals
\[N \subset N \Lambda'' \subset N\Lambda''N^{-1} = \Lambda'.\]
Since $N$ is a maximal left Ideal in $\Lambda'$ it follows that $N\Lambda''$ is either $N$ or $\Lambda'$.  If $N \Lambda'' = N$ then $\Lambda'' \subset O_r(N)$, contradiction.  So $N\Lambda'' = \Lambda'$ and thus
\[\Lambda' = N\Lambda''N^{-1} = \Lambda' N^{-1},\]
so $N^{-1}=\Lambda'$ which is impossible.  There is no counterexample.  $\blok$\\

{\defi Let $M,N$ be any pair of full $R$-lattices in $A$.  We say that their product $MN$ is \textbf{proper} if $O_r(M)=O_l(N)$.  Likewise, a product $M_1\cdot\ldots\cdot M_k$ is proper if $O_r(M_i)=O_l(M_{i+1})$, $\forall 1 \leq i \leq k-1$.}\\

{\st Let $M$ be a left Ideal of the maximal $R$-order $\Lambda$ in $A$.  Suppose that the $\Lambda$-module $\Lambda/M$ has composition length $k$.  Then $M$ is expressible as a proper product of $k$ maximal left integral Ideals $M_1, \ldots, M_k$ such that
\[O_l(M)=O_l(M_1), \ \ \ O_r(M)=O_r(m_k).\]}\\
\textbf{Proof} By induction on $k$.  Result is trivial when $k=1$, so assume $k > 1$, and suppose the theorem is true for the case $k-1$.  We can choose a maximal left Ideal $N$ of $\Lambda$ such that $M \subset N \subset \Lambda$, and then $N/M$ has composition length $k-1$ as a left $\Lambda$-module.  Both $M$ and $N$ are normal Ideals with left order $\Lambda$.  Using Theorem \ref{max13} we have
\[NN^{-1}=\Lambda, \ \ \ N^{-1}N =O_r(N) = \Lambda'.\]
There is a bijection $W \mapsto N^{-1}W$ between the set of left $\Lambda$-submodules $W$ of $N$ and left $\Lambda'$-submodules $N^{-1}W$ of $\Lambda'$, with inverse given by $N^{-1}W \mapsto N(N^{-1}W)$.  Therefore $N^{-1}M$ is a left Ideal of $\Lambda'$ and the quotient $\Lambda'/N^{-1}M$ has compostion length $k-1$ as a left $\Lambda'$-module.  Using the induction hypothesis we may write
\[N^{-1}M = M_2\cdot\ldots\cdot M_k\]
as a proper product of $k-1$ maximal integral Ideals, with
\[O_l(M_2)= O_l(N^{-1}M) = \Lambda', \ \ \ O_r(M_k)=O_r(N^{-1}M).\]
But then $M = N M_2 \cdot\ldots\cdot M_k$ is the desired proper product of $k$ maximal integral Ideals. $\blok$\\

\newpage
{\st \label{max16} With notations as in Lemma \ref{max6}
\begin{enumerate}
	\item A proper product of integral Ideals is integral.
	\item Let $M_{12}$ and $N_{34}$ be normal Ideals.  Then $M_{12} \subset N_{34}$ if and only if $M_{12}$ is a proper product
	\[M_{12}=B_{13}N_{34}C_{42},\]
	for some integral Ideals $B_{13}$ and $C_{42}$.
	\item Let $M_{12}$ and $N_{34}$ be normal Ideals with the same left order $\Lambda_1$.  Then $M_{12} \subset N_{14}$ if and only if $M_{12} = N_{14}C_{42}$ for some integral Ideal $C_{24}$.
	
\end{enumerate}
}
\begin{flushleft}\textbf{Proof}\end{flushleft}
\begin{enumerate}
	\item We show by induction on $k$ that 
	\[M_{12}M_{23}\cdot\ldots\cdot M_{k, k+1} \subset \Lambda_{k+1},\]
	if the $M$'s are integral Ideals.  This is clear for $k=1$ so let $k>1$ and assume that the result holds for $k-1$ factors.  So we find
	\[M_{12}\cdot\ldots\cdot M_{k-1,k}M_{k,k+1} \subset \Lambda_k M_{k,k+1} \subset \Lambda_{k+1}.\]
	\item First suppose that we have Ideals $B_{13}$ and $C_{42}$ such that
	\[M_{12}=B_{13}N_{34}C_{42},\]
	then we see
	\[M_{12}\subset \Lambda_3 N_{34} \Lambda_4 \subset N_{34}.\]
	Conversely, suppose $M_{12}\subset N_{34}$, and set
	\[B_{13}=M_{12}(\Lambda_3 M_{12})^{-1},\]
	\[C_{42} = (N_{34})^{-1}M_{12}.\]
	These Ideals are normal and have the indicated left and right orders.  We see that the normal Ideals $M_{12}$ and $\Lambda_3 M_{12}$ have the same right order, and $M \subset \Lambda_3 M$, so we obtain
	\[(M_{12})^{-1}\supset(\Lambda_3 M_{12})^{-1}.\]
	So
	\[B_{13} \subset M_{12}(M_{12})^{-1}=\Lambda_1,\]
	and this means $B_{13}$ is integral.  Further
	\[C_{42} = N_{34}^{-1}M_{12} \subset N_{34}^{-1}N_{34} \subset \Lambda_4,\]
	so $C_{42}$ is an integral Ideal.
	\begin{eqnarray*}
	B_{13}N_{34}C_{42}&=&M_{12}(\Lambda_3 M_{12})^{-1}\cdot N_{34}\cdot N_{34}^{-1}M_{12}\\
	&=&M_{12}(\Lambda_3 M_{12})^{-1}\cdot \Lambda_3 M_{12}\\
	&=& M_{12}\Lambda_2=M_{12}.
	\end{eqnarray*}
	\item If we take $\Lambda_1 = \Lambda_3$ in 2. we obtain $B=M_{12}(\Lambda_1 M_{12})=\Lambda_1$. $\blok$
\end{enumerate}
	
{\lem \label{max17} Let $M,N$ be normal Ideals in $A$.  Then $M \cdot N$ is a proper product if and only if replacing either factor by a larger $R$-lattice in $A$ increases the product.}\\
\textbf{Proof} In the notation of Lemma \ref{max6}, write $M=M_{12}$ and $N=N_{34}$.  Suppose that increasing either factor increases the product MN.  Then
\[MN = M\Lambda_2 N \Rightarrow N = \Lambda_2 N,\]
since we increase the factor but the product does not change.  Therefore $\Lambda_3 = \Lambda_2$.  Conversely, suppose $\Lambda_2 = \Lambda_3$, and suppose $M'N = MN$ where $M'$ is an $R$-lattice properly containing $M$.  Then
\[M' \subset M'\Lambda_2 = M' N N^{-1} = MNN^{-1} = M,\]
which is impossible. $\blok$	\\

{\st Let $M$ be an Ideal with $O_l(M)$ maximal order in $A$.  Then $M$ is a maximal left Ideal in $O_l(M)$ if and only if $M$ is a maximal right Ideal in $O_r(M)$.}\\
\textbf{Proof} We use notation as in Lemma \ref{max6}.  Let $M = M_{12}$ be a maximal left Ideal in $O_l(M) = \Lambda_1$.  Then $M \subset O_r(M)=\Lambda_2$ by Lemma \ref{max6}.  If $M = O_r(M)$, then $O_l(M)=O_r(M)=M$ which is impossible.  Suppose that $M$ is not a maximal right Ideal of $O_r(M)$.  Then there exists a normal Ideal $N=N_{32}$ such that
\[M \subset N \subset O_r(M),\]
where all inclusions are strict.  If $N = \Lambda_3$, then $\Lambda_2 = \Lambda_3 = N$, so we have a strict inclusion $N \subset \Lambda_3$.  From Theorem \ref{max16} we may write $M_{12} = B_{13} N$ for some integral Ideal $B=B_{13}$.  The product $B N$ is a proper product, so we find by Lemma \ref{max17}:
\[M = BN \subset B \Lambda_3 = B \subset \Lambda_1,\]
and $BN \neq B\Lambda_3$.  Then, since $M$ is maximal, $B = \Lambda_1$.  So $\Lambda_3 = O_r(B)=\Lambda_1$ which implies that $M = BN = \Lambda_3 N_{32} = N$.  Contradiction. $\blok$\\

{\st Let $\Lambda_1, \Lambda_2$ be a pair of maximal $R$-orders in $A$.  Then there exists a normal Ideal $M = M_{12}$.  If $I(\Lambda_j)$ the group of two-sided $\Lambda_j$-Ideals $j =1,2$ in $A$, there is an isomorphism
\[\phi_{12}:I(\Lambda_1) \rightarrow I(\Lambda_2): X \mapsto M^{-1}XM.\]
The map $\phi_{12}$ is independent of the choice of $M_{12}$.}\\
\textbf{Proof} The product $M = \Lambda_1 \cdot \Lambda_2$ is a normal Ideal with left order $\Lambda_1$ and right order $\Lambda_2$.  It is clear that $\phi_{12}$ carries $I(\Lambda_1)$ onto $I(\Lambda_2)$.  If $N_{12}$ is another normal Ideal, then $MN^{-1}\in I(\Lambda_1)$.  Since $I(\Lambda_1)$ is abelian, it follows that $MN^{-1}$ commutes with each $X \in I(\Lambda_1)$, and hence $M^{-1}XM = N^{-1}XM$.  Therefore $\phi_{12}$ does not depend on the choice of the normal Ideal $M_{12}$ used to define $\phi_{12}$.  Finally, $\phi_{12}$ is an isomorphism since it has an inverse
\[\phi_{21}: I(\Lambda_2)\rightarrow I(\Lambda_1):Y \mapsto MYM^{-1}.\]
This completes the proof. $\blok$ \\

\section{Graded Theory}

\subsection{Setting}

For this section, we will use the following notation and conventions.
\begin{itemize}
	\item $R$ is a Noetherian domain.
	\item $K = \Q(R)$ the quotient field of $R$.
	\item $G$ is a finite group.
	\item There is a group action of $G$ on $R$, denoted by $\sigma_g \in \textup{Aut}R$, where $g \in G$.
	\item $A$ is crystalline graded, $A = K \CGR G$.
\end{itemize}

Remark that the generating set $X_A$ (as used in Section \ref{max1}) in this case is $\{u_g|g\in G\}$.

\subsection{Lattices}

The definitions concerning lattices stay the same as in Section \ref{max1}.  We have the following however:

{\defi \label{max18} Let $T$ be any subset of $A$.  We call $T$ \textbf{graded} (in $A$) if and only if for each $x = \sum_{g \in G}x_g u_g \in T$, $x_g u_g \in T$.}\\

{\prop \label{max19} Let $L$ be a graded full $R$-lattice in $A$.  Then $L^{-1}$ as defined in Definition \ref{max3}, i.e.
\[L^{-1}=\{x \in A | LxL \subset L\}\]
is graded.}\\
\textbf{Proof} Let $x \in L^{-1}$, $x = \sum x_g u_g$.  Then for each $a = \sum a_g u_g$ and $b = b_g u_g \in L$ we have
\begin{eqnarray*}
a_g u_g x b_h u_h \in L \ \ \forall g,h \in G &\Rightarrow& a_g u_g \left(\sum_{t \in G} x_t u_t \right) b_h u_h \in L\ \ \forall g,h \in G\\
\ &\Rightarrow& a_g u_g x_t u_t b_h u_h \in L \ \ \forall g, h, t \in G\\
\ &\Rightarrow& \sum_{g \in G} \sum_{h \in G} a_g u_g x_t u_t b_h u_h \in L \ \ \forall t \in G\\
\ &\Rightarrow& a (x_t u_t) b \in L \ \ \forall t \in G,
\end{eqnarray*} 
and this means $L^{-1}$ is graded. $\blok$\\

\subsection{Orders}
Again we can use all definitions from Section \ref{max4}, with the added condition of Definition \ref{max18}:

{\defi We call a graded full $R$-latice $\Lambda$ in $A$, a \textbf{graded order or gr-order in $A$} if it is a graded ring with the same identity element as $A$.  A graded full $R$-lattice that is maximal with this property is called a gr-maximal order.}\\

{\prop If $M$ is a graded full $R$-order in $A$, $O_l(M)$ and $O_r(M)$ are graded orders.}\\
\textbf{Proof} That they are orders follows from Proposition \ref{max5}.  The fact that they are graded follows from a similar argument as in Proposition \ref{max19}. $\blok$\\

{\opm To make a distinction between a maximal object with the additional property of being graded (universally maximal), and a graded object that is maximal with the properties of such an object (maximal in all graded objects), we call the last one a gr-maximal.  This distinction makes sense, since there are gr-maximal orders which are not maximal and graded.  See \cite{FVO} for an example where a group ring is a gr-maximal order but not a maximal graded order.}\\

\subsection{Ideals}

{\defi Call a graded Ideal $M$ a \textbf{graded normal Ideal} if $O_l(M)$ is a gr-maximal $R$-order.  An \textbf{integral graded Ideal} is a graded normal Ideal such that $M \subset O_l(M)$.  A \textbf{gr-maximal integral Ideal} is a graded integral Ideal which is a gr-maximal left Ideal in $O_l(M)$.}\\

{\prop Suppose $I$ is a graded ideal (not necessarily a full lattice) of $A$.  Then $I$ is an Ideal of $A$, in other words, it is a full lattice.}\\
\textbf{Proof}  Let $a_g u_g \in I$ for some $g \in G$.  Then
\[\alpha^{-1}\left(g,g^{-1}\right)a_g u_g = a_g \alpha^{-1}\left(g,g^{-1}\right) u_g\in KI,\]
and $u_{g^{-1}} a_g \alpha^{-1}\left(g,g^{-1}\right) u_g\in AKI = KAI=KI$.  This implies that \[\sigma_{g^{-1}}(a_g)u_{g^{-1}}\alpha^{-1}\left(g,g^{-1}\right) u_g = \sigma_{g^{-1}}(a_g) \in KI.\]
Now $A \subset AKI = KI$.  $\blok$\\

{\defi \label{max20} A \textbf{gr-prime Ideal} of an $R$-order $\Lambda$ is a proper two-sided graded Ideal $P$ in $\Lambda$ such that for each pair of two-sided graded Ideals $S$ and $T$ in $\Lambda$ with $ST \subset P$ we have $S \subset P$ or $T\subset P$.}\\

{\prop \label{max21}With notations as in Definition \ref{max20} we can assume both $S$ and $T$ contain $P$ and as such $J\cdot J' = 0$ in $\Lambda/P$, where both $J$ and $J'$ are graded, then $J = 0$ or $J' =0$.}\\
\textbf{Proof} Suppose $S$ and $T$ fulfill Definition \ref{max20}.  Then $(S + P)(T+P) \subset ST + P \subset P$ and thus $S+P \subset P$ or $T+P \subset P$ and this means $S \subset P$ or $T \subset P$.  The other statement is now easy to prove.  $\blok$\\

For the rest of this section, we suppose $R$ is a Dedekind domain.

{\st The gr-prime Ideals in $\Lambda$ coincide with the gr-maximal two-sided Ideals of $\Lambda$.  If $P$ is a gr-prime Ideal of $\Lambda$ then $p = P \cap R^G = P_e$ is a nonzero prime in $R^G$ and $\bar \Lambda = \Lambda/P$ is a finite dimensional gr-simple $R^G/p$-ring.}\\
\textbf{Proof}  Let $P$ be a gr-prime Ideal in $\Lambda$ and set $p = R^G \cap P$, $\bar \Lambda = \Lambda/P$.  $P$ is a full $R$-lattice in $\Lambda$, so $p \neq 0$.  $p$ is proper since $1 \notin P$.  Now if $\alpha, \beta \in R^G$ then $\alpha \Lambda$ and $\beta \Lambda$ are graded two-sided Ideals of $\Lambda$.  Then
\begin{eqnarray*}
\alpha\beta \in p &\Rightarrow& \alpha\Lambda\beta\Lambda \subset P \Rightarrow \alpha\Lambda \subset P \textup{ or } \beta\Lambda \subset P\\
\ &\Rightarrow& \alpha \in p \textup{ or } \beta \in p,
\end{eqnarray*}
so $p$ is prime and $\bar\Lambda$ is finite dimensional over the field $R^G/p$.  So $\bar\Lambda$ is artinian, so it is gr-artinian, and thus the graded radical $\textup{rad}_{gr}\bar \Lambda$ is nilpotent.  This means by Proposition \ref{max21} that $\textup{rad}_{gr}\bar \Lambda = 0$.  In other words, $\bar \Lambda$ is gr-semisimple.  The gr-simple components are two-sided graded Ideals that annihilate each other, so again by Proposition \ref{max21} there is only one gr-simple component.  This means $\bar \Lambda$ is gr-simple and $P$ is gr-maximal.  The other way is trivial. $\blok$\\

{\lem Let $M$ be a graded Ideal in in a gr-maximal order $\Lambda$ of $A$.  Then $O_l(M)=O_r(M)=\Lambda$.  This implies that $M^{-1}$ is a graded $\Lambda$-Ideal in $A$.  Furthermore, if $N$ is another graded Ideal in $\Lambda$ with $N \subset M$, then $M^{-1} \subset N^{-1}$.}\\
\textbf{Proof} The proof is identical to that of Lemma \ref{max11}, while noting that $O_r(M)$ and $O_l(M)$ are graded. $\blok$\\

{\lem Let $M$ be a two-sided graded Ideal in the gr-maximal order $\Lambda$.  If $M \neq \Lambda$ then $M^{-1} \neq \Lambda$.}\\
\textbf{Proof} The proof is similar to that of Lemma \ref{max12}.  We note that in a gr-noetherian ring every graded Ideal contains a finite product of gr-prime Ideals. $\blok$\\

The following theorems, propositions and lemma's use nearly the same proofs as their nongraded counterparts.  We state them here without proof.\\

{\st Let $M$ be a full graded $R$-lattice in $A$ such that $O_l(M)$ is a gr-maximal order.  Then
\[M \cdot M^{-1} = O_l(M).\]}\\

{\st Let $\Lambda$ be a gr-maximal order.  For each two-sided graded Ideal $M$ in $\Lambda$ we have
\[MM^{-1}=M^{-1}M = \Lambda,\]
\[(M^{-1})^{-1}=M.\]
Furthermore, every two-sided graded Ideal in $\Lambda$ is uniquely expressible as a product of gr-prime Ideals of $\Lambda$, and multiplication of gr-prime Ideals is commutative.}\\

{\gev If $M$ is a proper left graded Ideal of the gr-maximal order $\Lambda$, then $M^{-1} \neq \Lambda$.}\\

{\st Let $M$ be a full graded $R$-lattice in $A$.  Then $O_l(M)$ is a gr-maximal order if and only if $O_r(M)$ is a gr-maximal order.}\\

{\defi Let $M,N$ be any pair of full graded $R$-lattices in $A$.  We say that their product $MN$ is \textbf{proper} if $O_r(M)=O_l(N)$.  Likewise, a product $M_1\cdot\ldots\cdot M_k$ is proper if $O_r(M_i)=O_l(M_{i+1})$, $\forall 1 \leq i \leq k-1$.}\\

{\st Let $M$ be a graded left Ideal of the gr-maximal $R$-order $\Lambda$ in $A$.  Suppose that the $\Lambda$-module $\Lambda/M$ has composition length $k$.  Then $M$ is expressible as a proper product of $k$ gr-maximal left integral Ideals $M_1, \ldots, M_k$ such that
\[O_l(M)=O_l(M_1), \ \ \ O_r(M)=O_r(m_k).\]}\\

{\st With notations as in Lemma \ref{max6}
\begin{enumerate}
	\item A proper product of graded integral Ideals is graded integral.
	\item Let $M_{12}$ and $N_{34}$ be graded normal Ideals.  Then $M_{12} \subset N_{34}$ if and only if $M_{12}$ is a proper product
	\[M_{12}=B_{13}N_{34}C_{42},\]
	for some graded integral Ideals $B_{13}$ and $C_{42}$.
	\item Let $M_{12}$ and $N_{34}$ be graded normal Ideals with the same graded left order $\Lambda_1$.  Then $M_{12} \subset N_{14}$ if and only if $M_{12} = N_{14}C_{42}$ for some graded integral Ideal $C_{24}$.
	
\end{enumerate}
}
	
{\lem Let $M,N$ be graded normal Ideals in $A$.  Then $M \cdot N$ is a proper product if and only if replacing either factor by a larger graded $R$-lattice in $A$ increases the product.}\\

{\st Let $M$ be a graded Ideal with $O_l(M)$ a gr-maximal order in $A$.  Then $M$ is a gr-maximal left Ideal in $O_l(M)$ if and only if $M$ is a gr-maximal right Ideal in $O_r(M)$.}\\

{\st Let $\Lambda_1, \Lambda_2$ be a pair of gr-maximal $R$-orders in $A$.  Then there exists a graded normal Ideal $M = M_{12}$.  If $I(\Lambda_j)$ the group of two-sided graded $\Lambda_j$-Ideals $j =1,2$ in $A$, there is an isomorphism
\[\phi_{12}:I(\Lambda_1) \rightarrow I(\Lambda_2): X \mapsto M^{-1}XM.\]
The map $\phi_{12}$ is independent of the choice of $M_{12}$.}\\

\end{document}